\documentclass[a4 paper, 12pt]{article}

\usepackage{amsfonts}
\usepackage{amssymb}
\usepackage{amsmath}
\usepackage{fancyhdr}
\usepackage{eufrak}
\usepackage{latexsym}
\usepackage{amsbsy}
\usepackage[all]{xy}
\usepackage[latin1]{inputenc}
\usepackage[T1]{fontenc}
\usepackage{amsmath}

\newenvironment{theorem}[1]
{\vskip 2mm\noindent \textbf{Theorem #1}  \it}{\vskip 2mm}

\newenvironment{definition}[1]
{\vskip 2mm\noindent \textbf{Definition #1}  \it}{\vskip 2mm}

\newtheorem{defn}{Definition}[section]
\newtheorem{thm}[defn]{Theorem}
\newtheorem{prop}[defn]{Proposition}
\newtheorem{lem}[defn]{Lemma}

\newtheorem{rem}[defn]{Remark}

\newcommand{\proof}{\vskip 2mm \noindent {\textsc{Proof: }}\rm}
\newcommand{\proo}{\vskip 2mm \noindent {\textsc{Proof }}\rm}
\newcommand{\fin}{\hfill{\Large$\Box$}\\}
\newcommand{\finsec}{\hfill{\Large$\Box$}}

\newcommand{\al}{\alpha}

\newcommand{\Ga}{\Gamma}
\newcommand{\ga}{\gamma}
\newcommand{\ale}{\textrm{-a.e.}}
\newcommand{\si}{\sigma}

\newcommand{\om}{\omega}\newcommand{\epsi}{\epsilon}

\newcommand{\C}{\mathbb {C}}
\newcommand{\R}{\mathbb {R}}
\newcommand{\N}{\mathbb {N}}
\newcommand{\Z}{\mathbb {Z}}
\newcommand{\E}{\mathbb {E}}

\newcommand{\Pj}{\mathbb {P}}
\newcommand{\cdb}{\textsf {Card }}

\newcommand{\Id}{{\rm Id}}
\newcommand{\Diam}{{\textsf {diam} \, }}
\newcommand{\Leb}{{\textsf {Leb}  }}

\newcommand{\Lebe}{{\textsf {Leb}\,  }}
\newcommand{\Diag}{{\textsf {Diag}  }}
\newcommand{\sph}{{\textsf s }}
\newcommand{\area}{{\textsf {area}  }}
\newcommand{\per}{{\textsf {Per}  }}

\newcommand{\Vol}{{\rm Vol \, }}

\newcommand{\Jac}{{\textsf {Jac} \,  }}
\newcommand{\conf}{{\textsf {conf} \,  }}

\newcommand{\degr}{{\rm deg \, }}
\newcommand{\Mod}{{\textsf {mod}  \, }}

\newcommand{\tal}{{\tilde \alpha}}
\newcommand{\tbe}{{\tilde \beta}}

\def\abs#1{\vert #1\vert}
\def\norm#1{\left\|\, #1\,\right\|}

\def\TT{{\cal T}}

\def\EE{{\cal E}}

\def\AA{{\cal A}}

\def\BB{{\cal B}}
\def\CC{{\cal C}}

\def\DD{{\cal D}}

\def\FF{{\cal F}}
\def\UU{{\cal U}}
\def\GG{{\cal G}}
\def\JJ{{\cal J}}
\def\RR{{\cal R}}
\def\II{{\cal I}}
\def\LL{{\cal L}}
\def\VV{{\cal V}}
\def\MM{{\cal M}}
\def\GG{{\cal G}}
\def\SS{{\cal S}}
\def\TT{{\cal T}}

\def\QQ{{\cal Q}}
\def\NN{{\cal N}}
\def\WW{{\cal W}}

\def\HH{{\cal H}}
\def\com{\ar@{}[rd]|{\circlearrowleft}}

\title {Bernoulli coding map and almost sure invariance principle for endomorphisms of $\Pj^k$}
\author{C. Dupont}
\date{ \today }

\begin{document}

\maketitle

\begin{abstract}
Let $f$ be an holomorphic endomorphism of $\mathbb{P}^k$ and $\mu$ be its measure
of maximal entropy. We prove an Almost Sure Invariance Principle for
the systems $(\mathbb{P}^k,f,\mu)$. Our class $\UU$ of observables includes the
H\"older functions and unbounded ones which present
analytic singularities. The proof is based on a geometric construction of a Bernoulli coding map $\omega: (\Sigma , s , \nu) \to
(\mathbb{P}^k,f,\mu)$. We obtain the invariance principle for an observable $\psi$ on
$(\mathbb{P}^k,f,\mu)$ by applying Philipp-Stout's theorem for $\chi = \psi
\circ \omega$ on $(\Sigma , s , \nu)$. The invariance principle implies the Central Limit Theorem as well as several statistical properties for the class $\UU$. As an application, we give a
\emph{direct} proof of the absolute continuity of the measure $\mu$
when it satisfies Pesin's formula. This approach relies on the
Central Limit Theorem for the unbounded observable $\log
\textsf{Jac} \,  f \in \UU$.
\end{abstract}

\small{\noindent \emph{Key Words:} holomorphic dynamics, Bernoulli coding map, almost sure invariance principle.

}\normalsize

\section{Introduction}\label{intro}

Let $f : \Pj^k \to \Pj^k$ be an holomorphic endomorphism of algebraic degree $d
\geq 2$. Its equilibrium measure $\mu$ is the limit of the probability measures $d_t^{-n} (f^n)^* \eta^k$,
where $d_t:= d^k$ is the topological degree of $f$ and $\eta^k$ is the
standard volume form on $\Pj^k$. We refer to the survey article of
Sibony \cite{Sib} for an introduction to the dynamical systems $(\Pj^k,f,\mu)$. Fornaess-Sibony proved that $\mu$
is mixing \cite{FS} and Briend-Duval that $\mu$ is the unique measure of maximal entropy \cite{BrDu2}.  \\

 Przytycki-Urba\'nski-Zdunik \cite{PUZ} introduced coding techniques
 for $(\Pj^1,f,\mu)$. This allowed them to prove the Almost Sure
 Invariance Principle (ASIP) for H\"older and
 singular observables, like $\log \abs{f'}$. In the present article, we
 extend the coding techniques to $(\Pj^k,f,\mu)$ and obtain the
 ASIP for observables which allow analytic singularities. As an application, we
 obtain a direct proof of the absolute continuity of $\mu$
 when it satisfies Pesin's formula. We review our results in subsections \ref{baba} - \ref{jaja}, subsection \ref{relt} is devoted to related results.

\subsection{Bernoulli coding maps}\label{baba}

Let us endow $\Sigma
:=  \{ 1, \ldots, d_t \}^\N$ with the natural product measure $\nu:= \otimes_{n = 0}^\infty \bar
\nu$, where $\bar \nu$ is equidistributed on $\{ 1, \ldots, d_t \}$. We
 denote by $\tal$ the elements of
$\Sigma$ and by $s$ the left shift acting on $\Sigma$. Let $\JJ$ be the support
 of $\mu$. The following theorem yields coding maps $\om
: \Sigma \to \JJ$ up to zero measure sets. The set $\SS \subset \Pj^k$ will be
 defined in section \ref{GCT}, it has zero Lebesgue measure.
\begin{theorem}{A:}
Let $z \in \Pj^k \setminus \SS$. There exist an $s$-invariant set $\Sigma'
\subset \Sigma$ of full $\nu$-measure and an $f$-invariant set $\JJ'
\subset \JJ$ of full $\mu$-measure satisfying the following properties. For any $\tal \in
\Sigma'$, the point $\om(\tal):= \lim_{n \to \infty} z_n(\tal) \in \JJ'$ is
well defined. We have $\om_* \nu  = \mu$ and the following diagram commutes:
\[ \xymatrix{
           \Sigma'  \ar[rr]^{s}  \ar[d]_\om &    & \Sigma'  \ar[d]^\om \\
           \JJ'  \ar[rr]^{f}       & &  \JJ'        }
\]
Moreover there exist $\theta , \epsi > 0$, $n_z \geq 1$ and $\tilde n : \Sigma' \to \N$ larger than $n_z$ such that:
\begin{enumerate}
\item $d(z_n(\tal),\om(\tal) ) \leq \tilde c_\epsi \, d^{-\epsi n}$  for every $\tal \in \Sigma'$ and $n \geq \tilde n (\tal)$,
\item $\nu ( \{ \tilde n  \leq q \} )  \geq 1 - c_\theta \, d^{-\theta q}$ for every $q \geq n_z$.
\end{enumerate}
\end{theorem}

We note that $\Sigma'$, $\JJ'$ and $\om$
depend on $z \in \Pj^k \setminus \SS$, but $\om_* \nu = \mu$ holds
true for any such $z$. Observe also that $\om$ is not necessarily
injective. The proof of theorem A (see section \ref{GCT}) is based on the construction of a geometric coding
tree, following the approach of Przytycki-Urba\'nski-Zdunik \cite{PUZ}
for $(\Pj^1,f,\mu)$. The point $z$ is the \emph{root} of the tree, and the set $\{ z_n(\tal) , \tal \in \Sigma \}$ is a suitable enumeration of the
$d_t^{n+1}$ points of $f^{-(n+1)}(z)$, these are \emph{vertices} of the
tree. The convergence of $(z_n(\tal))_n$ for a generic $\tal \in \Sigma$ is obtained by constructing $d_t$ good paths joining
$z$ to $w \in f^{-1}(z)$, whose inverse images decrease exponentially. In the context of
$(\Pj^1,f,\mu)$, that property was obtained in \cite{PUZ} by using Koebe
distortion theorem. The difficulty  in higher
dimensions is
to substitute this argument. We
establish for that purpose a quantified version of a theorem of
Briend-Duval (see section \ref{BD}).

\subsection{The class $\UU$ and approximation by cylinders}\label{byby}

\begin{definition}{:} An observable $\psi : \Pj^k \to \R \cup \{ -
  \infty \}$ belongs to the class $\UU$ if:

- $e^\psi$ is $h$-H\"older for some $h > 0$,

- $\NN_\psi:= \{ \psi = -\infty \}$ is a (possibly empty) proper algebraic set of
  $\Pj^k$,

- $\psi \geq \log d(\cdot,\NN_\psi)^\rho$ for some $\rho >0$.
\end{definition}
\noindent For instance, the H\"older functions are in $\UU$, as well
as the unbounded
function $\log \Jac f$. We will show that $\UU \subset
L^p(\mu)$ for any $1 \leq p < + \infty$ (see subsection \ref{cccd}).

\begin{theorem}{B:} Let $\psi \in \UU$ be a
  $\mu$-centered observable and $\om$ be a coding map provided by theorem A. Let $\chi:= \psi \circ
  \om$ and  $1 \leq p < + \infty$. We denote by $\E ( \chi \vert
  \CC_n  )$ the conditional expectation of $\chi$ with respect
  to the $(n+1)$-cylinders.
\begin{enumerate}
\item there
  exist $\hat c_p, \lambda_p > 0$ such that $\norm { \chi   - \E ( \chi \vert
  \CC_n  )  }_p \leq \hat c_p \, e^{-n\lambda_p}$ for every $n \geq 0$.

\item $R_j(\chi):= \int_\Sigma \chi \cdot
\chi \circ s^j \, d\nu$ satisfies $\abs{R_j(\chi)} \leq 2 \norm {\chi}_2  \hat c_2 \,
  e^{-(j-1) \lambda_2}$ for every $j \geq 1$.
\end{enumerate}
\end{theorem}
The proof occupies section \ref{estiho}, it is based on the
  regularity properties of $\om$ (namely the points 1, 2 of theorem A) and on the fact that $\mu$ is a Monge-Amp\`ere mass with H\"older
potentials. Theorem B allows us to prove theorem C below.

\subsection{Almost Sure Invariance Principle}\label{bibi}

Let $\psi \in L^2(\mu)$ be a $\mu$-centered observable and $S_n(\psi):= \sum_{j=0}^{n-1} \psi \circ f^j$. We say that $\psi$ satisfies the Almost Sure Invariance Principle (ASIP) if there exist, on an extended probability space, a
sequence of random variables $(\SS_n)_{n \geq 0}$ together with a
Brownian motion $\WW$ such that for some $\gamma > 0$:
\begin{enumerate}
\item[-]  $\SS_n =  \WW(n) + o(n^{1/2 -
  \gamma})$ almost everywhere,
\item[-] $(S_0(\psi),\ldots,S_n(\psi))$ and $(\SS_0,\ldots,\SS_n)$ have the
same distribution for any $n \geq 0$.
\end{enumerate}
\noindent We shall denote $\si$-ASIP to specify the variance of Brownian motion.

\begin{theorem}{C:} For every $\mu$-centered observable $\psi \in \UU$, we have:
\begin{enumerate}
\item $\si:= \lim_{n \to \infty} {1 \over \sqrt n} \norm{S_n(\psi)}_2$ exists, and $\si^2 = \int_{\Pj^k} \psi^2 \, d\mu + 2\sum_{j \geq 1} \int_{\Pj^k}
 \psi \cdot \psi \circ f^j \, d\mu$.

\item If $\si = 0$, then $\psi = u - u \circ f$ holds $\mu$-a.e. for some $u \in L^2(\mu)$.

\item If $\si > 0$, then $\psi$ satisfies the $\si$-ASIP.
\end{enumerate}
\end{theorem}

The ASIP implies classical limit theorems related to Brownian
 motion: the
 Central Limit Theorem (CLT), the Law of Iterated Logarithm, Kolmogorov integral tests (see \cite{De}, \cite{PS}). The ASIP also
implies the almost sure version of the CLT, meaning that ${1 \over
 \log n} \sum_{k=1}^n {1 \over k} \delta_{ {1 \over
    \sqrt{k}} S_k(\psi)(x)}$ converges $\mu \ale$ to the normal law
 $\NN(0,\si^2)$  (see \cite{LP}, \cite{CG}). \\

Let us outline the proof
 of theorem C (see section \ref{CC}). Let $\om : \Sigma \to \Pj^k$ be a coding map provided by theorem A and $\psi \in \UU$. Since $\om$ satisfies $f
\circ \om = \om \circ s$ and $\om_*\nu=\mu$, we are reduced to prove
the assertions for $\chi = \psi \circ \om$ on $(\Sigma,s,\nu)$. The points 1 and 2 follow from theorem
 B(2) and classical arguments. The point 3 is a consequence of theorem
 B(1) and Philipp-Stout's theorem (\cite{PS}, section 7). That result relies on an approximation of the partial sums of $(\chi
\circ s^j)_{j \geq 0}$ by a sequence of martingale differences defined with respect to the increasing
filtration $(\CC_n)_{n \geq 0}$.

\subsection{An application to smooth ergodic theory}\label{jaja}

Let  $\chi_1 \leq \ldots \leq \chi_k$ be the Lyapunov exponents  of
$\mu$. Briend-Duval  \cite{BrDu1} proved that they are larger than or equal to
$\log d^{1/2}$. Since $\mu$ has entropy $\log d^k$, Pesin's formula $h(\mu) = 2 (\chi_1 + \ldots + \chi_k)$
holds if and only if these exponents are minimal. We proved in a previous
article that $\mu$ is then absolutely
continuous with respect to Lebesgue measure \cite{Du}. We there followed
the classical approach of Sinai-Pesin-Ledrappier, based on
the construction of a suitable invariant partition which is dilated
and realizes entropy (see \cite{P1}, \cite{Le}). We propose in section \ref{appll} a new proof,
based on the CLT for the unbounded $\mu$-centered observable $J:= \log
\Jac f -  2 (\chi_1 + \ldots + \chi_k) \in \UU$. We obtain the following result,
where $\si_J: = \lim_{n \to \infty} {1 \over \sqrt n}
\norm{S_n(J)}_2$.
\begin{theorem}{D:}
If the Lyapunov exponents are minimal equal to $\log d^{1/2}$, then
$\si_J = 0$, and $\mu$ is absolutely continuous with respect to
Lebesgue measure.
\end{theorem}
A crucial fact for the proof is that for any holomorphic endomorphism of $\Pj^k$ and any $\mu$-generic point $x \in \Pj^k$, the
minimal dilation rate of $f^n$ at $x$ (i.e. $\norm{(d_x
  f^n)^{-1}}^{-1}$) is bounded below by $d^{n/2}$ up to the multiplicative factor $1/n$. In other words, the usual
$e^{-n\epsi}$-correction, due to the non-uniform hyperbolicity of
$(\Pj^k,f,\mu)$, can be replaced here by $1/n$. This was proved by
Berteloot-Dupont \cite{BeDu}, using a pluripotential result of
Briend-Duval \cite{BrDu1} and the fact that $\mu$ is a Monge-Amp\`ere
mass. In particular, the
product of the dilation rates satisfies $\Jac
f^n(x) \geq \norm{(d_x f^n)^{-1}}^{-2k}  \geq (d^{n/2}/n)^{2k} = d^{kn}/n^{2k}$. Now if we assume $\si_J > 0$, then the function $\log \Jac f^n$ would present non trivial oscillations
around its mean value $\log d^{kn}$, due to the CLT. More precisely, it would imply $\log \Jac f^n \leq \log d^{kn} - \si_J \sqrt n$  on a subset of $\mu$-measure $\simeq \int_{-\infty}^{-1}
e^{-u^2/2}$. That contradicts the preceding estimate, hence $\si_J  = 0$. We deduce the
absolute continuity of $\mu$ from the cocycle property $J = u - u
\circ f$ $\mu$-a.e. and a linearization property of the dynamics along typical
negative orbits \cite{BeDu}.

\subsection{Related results}\label{relt}

The systems $(\Pj^k , f , \mu)$ and $(\Sigma , s
, \nu)$ are actually conjugated by a bimeasurable map up to zero measure
subsets, that property was proved by Briend \cite{Br}. However, the regularity of the conjugacy seems difficult to
handle. Let us also mention that finite-to-one coding maps $(\Pj^k,f,\mu) \to (\Sigma,s,\nu)$ were constructed by Buzzi \cite{Bu} by means of suitable partitions of $\Pj^k$. \\

The ASIP has been proved for many dynamical systems: for piecewise monotonic maps by Hofbauer-Keller
\cite{HK}, for Anosov maps by Denker-Philipp \cite{DP} and for partially and
non-uniformly hyperbolic systems by Dolgopyat \cite{Do} and
Melbourne-Nicol \cite{MN}. We refer to the survey articles of
Chernov \cite{C} and Denker \cite{De} for limit theorems and
statistical properties concerning dynamical systems. \\

The ASIP implies the CLT. Nevertheless, the latter can be directly proved
 via coding techniques and Ibragimov's theorem \cite{I}. That
 method was employed by Sinai \cite{Sin}
and Ratner \cite{R} for the geodesic flow in negative
curvature, and by Bowen \cite{Bo} for Anosov maps. In the present article, Ibragimov's condition is fulfilled by theorem B. \\

The Gordin's theorem provides another method for proving the CLT (see \cite{G}, \cite{Li}). It relies on an approximation of $(\psi \circ
f^j)_{j \geq 0}$ by a sequence of reverse martingale differences. In our context, this can be done if $\sum_{n \geq 0} \norm{\Lambda^n
   \psi}_2$  (denoted $(\star)$) converges, where $\Lambda$ denotes the Ruelle-Perron-Frobenius
 operator (we have $\Lambda^n \psi (z) = {1 \over d_t^n}
 \sum_{y \in f^{-n}(z)} \psi(y)$ for every $z \in \Pj^k$). Let us note that the reverse martingale mentioned
 is defined with respect to the decreasing
 filtration $(f^{-n}\BB)_{n \geq 0}$, where $\BB$ is the Borel
$\si$-algebra of $\Pj^k$.

The exponential decay of
correlations ensures the convergence of $(\star)$. This was proved in the context of $(\Pj^k,f,\mu)$ by Fornaess-Sibony \cite{FS} for $C^2$
observables and by Dinh-Sibony for H\"older observables
\cite{DS2}. Dinh-Nguyen-Sibony have recently extended that property for differences of quasi-plurisubharmonic functions (the so-called \emph{dsh} functions) \cite{DNS2}. The proof relies on exponential estimates for plurisubharmonic functions with respect to $\mu$. They also obtained in that article a Large Deviations Theorem for bounded dsh and H\"older observables.  In \cite{DNS1}, Dinh-Nguyen-Sibony proved  the local CLT for $(\Pj^1,f,\mu)$ by using the theory of perturbed operators.

Denker-Przytycki-Urba\'nski \cite{DPU} employed a geometric method to prove the convergence of $(\star)$ for $(\Pj^1,f,\mu)$ and H\"older
observables. The idea was to compare $\Lambda^n \psi(z)$ to $\Lambda^n
\psi(z')$ by using the contraction of most of the inverse branches
of $f^n$. The cornerstone is a precise analysis of the dynamics
near the critical points in the support of $\mu$. Cantat-Leborgne \cite{CL}
 extended this approach to $(\Pj^k,f,\mu)$. A crucial ingredient was a polynomial estimate for the
 $\mu$-measure of postcritical neighbourhoods (lemma 5.7 of
 \cite{CL}). The original proof of that lemma contains a gap, the authors have recently proposed another one. Cantat-Leborgne also established in \cite{CL} a quantified version of the Briend-Duval theorem. Our version is similar, but we shall give a different proof.  \\

The systems $(\Pj^k,f,\mu)$ whose measure $\mu$ is absolutely
continuous with respect to Lebesgue measure were characterized by
Berteloot, Dupont and Loeb \cite{BeDu}, \cite{BL}. In that case, $f$
is semi-conjugated to an affine dilation on a complex torus, these
maps are the so-called \emph{Latt\`es examples}. We note that
theorem D characterizes these maps by the minimality of the Lyapunov
exponents. Another characterization of Latt\`es examples involves
the \emph{Hausdorff dimension} of $\mu$, defined as the infimum of
the Hausdorff dimension of Borel sets with full $\mu$-measure (see
Pesin's book \cite{P2}): Dinh-Dupont \cite{DD} proved that
$\dim_\HH(\mu) = 2k$ if and only if the exponents are minimal. In
the context of  $(\Pj^1,f,\mu)$, Ma\~n\'e \cite{Ma\~n} proved that
$\log d = \dim_\HH(\mu) \cdot \chi$, where $\chi$ denotes the
Lyapunov exponent of $\mu$. In particular, the function $L:= \log d
- \dim_\HH(\mu) \cdot \log \abs{f'}$ is a $\mu$-centered observable.
Zdunik \cite{Z} proved that $\si_L = 0$ if and only if $f$ is a
Latt\`es example, a Tchebychev polynomial or a power $z^{\pm d}$.
The proof relies on the classification of critically finite
fractions with parabolic Thurston's orbifold.

\section{Generalities}\label{cons}

\subsection{The holomorphic systems $(\Pj^k,f,\mu)$}\label{consa}

We introduce in this section the systems $(\Pj^k,f,\mu)$. We refer to the articles
\cite{BrDu1}, \cite{BrDu2}, \cite{FS} and \cite{Sib} for
definitions and properties. Here $\Pj^k$ denotes the complex projective space of dimension $k$. We denote by $\eta$ the Fubini-Study form on $\Pj^k$. This is a $(1,1)$-form defined in homogeneous coordinates by ${i \over 2 \pi} \partial  \bar \partial \log \norm{z}^2$. It induces the standard metric on $\Pj^k$, the volume of $\Pj^k$ with respect to this metric is equal to $1$. The form $\eta$ induces on every complex line $L \subset \Pj^k$ the spherical metric with area $1$. Let $f$ be an holomorphic endomorphism of
$\Pj^k$ with algebraic degree
$d \geq 2$. It is defined in homogeneous coordinates by $[ P_0 :
  \cdots : P_k ]$ where the $P_i$ are homogeneous
polynomials of degree $d$ (without common zero except the origin). The topological degree
of $f$ is $d_t:= d^k$. An inverse branch of
$f^n$ on $U \subset \Pj^k$ is an injective holomorphic map $g_n$
 satisfying $f^n \circ g_n = \Id_U$. We let $\per f:= \cup_{n \geq 1} \{ x \in \Pj^k
\, , \, f^n(x) = x \}$, this set is at most countable. Let $\CC$ be the
critical set of $f$, $\VV:= \cup_{i = 0}^\infty
f^i(\CC)$ and $\VV_n:= \cup_{i = 1}^ n f^i(\CC)$. The degree of
$\VV_n$, denoted $\tau_n$, is equal to $(d + \ldots +
d^n) \deg \CC$ counted with multiplicity.  \\

The equilibrium measure $\mu$ is defined as
the limit of $\mu_{n,z}:=
{1 \over d_t ^n} \sum_{f^n(y) = z} \delta_y$, where $\delta_y$ denotes
the Dirac mass at $y$. In that definition, $z$ has to be taken outside a totally invariant
algebraic set $\EE \subset \VV$, the so-called
exceptional set of $f$. We denote by $\JJ$ the support of $\mu$. The measure $\mu$ is mixing and satisfies $\mu(f(B)) = d_t
\, \mu(B)$ whenever $f$ is injective on $B$. It is the unique measure of
maximal entropy (equal to $\log d_t$). The Lyapunov exponents
$\chi_1 \leq \ldots \leq \chi_k$ of $\mu$ are larger than or equal to $\log
d^{1/2}$. They satisfy the classical formula $\int_{\Pj^k} \log \Jac f \, d\mu = 2 (
\chi_1 + \ldots +\chi_k )$, where $\Jac f$ is the non-negative $\CC^\infty$
function on $\Pj^k$ satisfying $f^* \eta^k = \Jac f
\cdot \eta^k$. The latter is the real jacobian of $f$, it vanishes on the critical set $\CC$ of $f$. \\

The measure
$\mu$ can also be defined \emph{via} pluripotential
theory: we have $\mu = T^k$, where $T$ is the Green current of
$f$. The latter is a closed positive
$(1,1)$ current on $\Pj^k$ with H\"older potentials. In particular, for any algebraic subset
$A \subset \Pj^k$, there exist $c,\ga > 0$ such that the
$r$-neighbourhood of $A$ satisfies $\mu ( A[r] )
\leq c \, r^\ga$ for any $r >
0$ (see \cite{DS4}, Prop. 2.3.7). For any
$\delta > 0$ and $\tilde c > 0$, we set $c_\delta:= (1-d^{-\delta})^{-1}$ and $\tilde c_\delta:= \tilde c (1-d^{-\delta})^{-1}$. In the sequel, $c > 0$ is a constant independent of $n$, it may differ from a line to another.

\subsection{The class $\UU$}\label{cccd}

Let us recall the definition of the class $\UU$ (see subsection \ref{byby}).
\begin{defn}\label{unb}
Let $\UU$ be the set of functions $\psi : \Pj^k \to \R \cup \{ -
\infty \}$ satisfying:

- $e^\psi$ is  $h$-H\"older on $\Pj^k$ for some $h >0$,

- $\NN_\psi:= \{ \psi = -\infty \}$ is a (possibly empty) proper algebraic set of
  $\Pj^k$,

- $\psi \geq \log d(\cdot,\NN_\psi)^\rho$  on $\Pj^k$ for some $\rho >0$.
\end{defn}
The H\"older functions belong to $\UU$. Examples of unbounded
observables are:
\begin{enumerate}
\item[-] the functions $\psi = \log \abs Q - q \log \norm{\cdot}$,
  where $Q$ is a $q$-homogeneous polynomial on $\C^{k+1}$. Here the algebraic subset $\NN_\psi$
  is the zero set of $Q$.
\item[-] the functions $\psi = \log \norm{\Lambda^j d_x f}$ ($1
  \leq j \leq k$), where $\Lambda^j d_x f$ is the $j$-exterior power of the differential $d_x
  f$. In particular, $\log \Jac f \in \UU$ (take $j = k$).
\end{enumerate}
The conditions of definition \ref{unb} are easy to verify
for these functions, the last one is a
consequence of Lojasiewicz's inequality (see \cite{Lo}, §4.7). We
prove below that $\psi \in L^p(\mu)$ for any $\psi \in \UU$ and
$1 \leq p < + \infty$. Actually, we establish an estimate for $\int_{\NN_\psi[r]}  \abs{ \psi }^p$, useful to prove
theorem B. We recall that
$\mu \left( \NN_\psi[r]  \right) \leq c \, r^\ga$ for some $c,\ga > 0$ (see subsection \ref{consa}).

\begin{prop}\label{cuun}
Let $\psi \in \UU$ and $1 \leq p < + \infty$. There exists $\kappa > 0$ such that:
\begin{equation*} \label{dec}
\forall \, 0 < r < 1/2 \ \ , \ \ \int_{\NN_\psi[r]}  \abs{ \psi }^p \, d\mu  \leq \kappa \, r^{\ga /2}.
\end{equation*}
In particular $\psi \in L^p(\mu)$.
\end{prop}

\proof Let $\psi \in \UU$ and $\NN:= \NN_\psi$. We may
assume that $0 \leq e^\psi \leq 1$ by adding some constant to $\psi$. Let $r < 1/2$ and $\QQ_j:= \NN[r/2^j] \setminus \NN[r/2^{j+1}]$. Since $e^\psi \geq (r/2^{j+1})^\rho$ on $\QQ_j$, we obtain:
\[ \int_{\NN[r]}  \abs {\psi }^p \, d\mu  = \sum_{j \geq 0}
\int_{\QQ_j}  \abs{ \log e^\psi} ^p \, d\mu \leq \sum_{j \geq 0}
\left \vert \, \rho \log \left( {r \over 2^{j+1}} \right) \, \right  \vert ^p \cdot \mu (\QQ_j).          \]
The inequalities $\mu (\QQ_j) \leq c (r/2^j)^\ga$ and $\abs{\log {r \over 2^{j+1}}} =  (j+1) \log 2 + \log {1 \over r}
 \leq   (j+2)   \log {1 \over r}$ yield:
\[ \int_{\NN[r]}  \abs{ \psi }^p \, d\mu \leq  \left[ c \, \rho^p \,  \sum_{j \geq
  0} {(j+2)^p \over 2^{\ga j}} \right]  \left( \log {1 \over r}
\right)^p   r^\ga = M_{\rho,\ga} \cdot \left(  \log {1 \over r} \right)^p  r^{\ga/2} \cdot  r^{\ga/2} .   \]
The lemma follows with $\kappa:=   M_{\rho,\ga} \cdot \sup_{ 0 < r < 1/2 } \left( \log {1 \over r} \right)^p r^{\ga/2}$. \finsec

\subsection{The Bernoulli space $(\Sigma,s,\nu)$}\label{consaaa}

We endow $\AA:= \{ 1, \ldots, d_t \}$ with the equidistributed
probability measure $\bar \nu$. We set $\Sigma:= \AA^\N$, $s : \Sigma \to \Sigma$ the left
shift and $\nu:= \otimes_{n =
  0}^\infty \bar \nu$. We denote by $\tal:= (\al_n)_{n \geq 0}$ the
elements of $\Sigma$, by $\CC_n$ the set of cylinders of
length $n+1$, and by $\pi_n : \Sigma \to \AA^{n+1}$ the projection $\pi_n
(\tal ):= (\al_0 , \ldots , \al_n )$. For any $\tal \in \Sigma$, we set
$C_n(\tal): = \pi_n^{-1} (\al_0 , \ldots , \al_n
)$. We denote by $\E(\chi \vert \CC_n)$ the
conditional expectation of $\chi \in L^2(\nu)$ with respect to $\CC_n$. If $\LL = \{ A_1 , \ldots
, A_p \} \subset \CC_n$, we set $\LL^*:= \cup_{1 \leq j \leq p} A_j$.

\subsection{Almost Sure Invariance Principle}\label{consb}

Let $(X,g,m)$ be either $(\Sigma,s,\nu)$ or $(\Pj^k,f,\mu)$. For any observable $\varphi \in
L^2(m)$, we set $S_n(\varphi):= \sum_{j = 0}^{n-1} \varphi \circ
g^j$ and $R_j(\varphi):= \int_X \varphi \cdot \varphi \circ g^j \,
dm$. We say that $\varphi$ is $m$-centered if $\int_X \varphi \, dm =
0$ and that $\varphi$ is a cocycle if $\varphi = u - u \circ g$ $m$-a.e. for some $u \in L^2(m)$. \\

An observable $\varphi$ on $(X,g,m)$ satisfies the Almost Sure Invariance Principle (ASIP) if there exist on a probability space
$(\tilde X , \tilde m)$ a sequence of random variables $(\SS_n)_{n \geq 0}$ and a Brownian motion $\WW$ such that:
\begin{enumerate}
\item[-]  $\SS_n =  \WW(n) + o(n^{1/2 -
  \gamma})$ $\tilde m$-a.e. for some $\gamma > 0$,
\item[-] $(S_0(\psi),\ldots,S_n(\psi))$ and $(\SS_0,\ldots,\SS_n)$ have the
same distribution for any $n \geq 0$.
\end{enumerate}
We denote $\si$-ASIP to specify the variance of Brownian
motion. The $\si$-ASIP implies the $\si$-Central Limit Theorem ($\si$-CLT),
meaning that:
\[ \forall t \in \R \ , \ \lim_{n \to \infty}  m \left( {S_n(\varphi)
  \over \si \sqrt n} \leq t  \right)  =  {1 \over \sqrt{2\pi}}
\int_{-\infty}^t e^{-u^2/2} \, du . \]

\begin{rem}{} Suppose that $\om : \Sigma \to \Pj^k$ is a coding map
  provided by theorem A. Since $\om_* \nu = \mu$ and $f \circ \om
= \om \circ s$, a $\mu$-centered observable $\psi \in L^{2}(\mu)$ satisfies the $\si$-ASIP if and only if the
$\nu$-centered observable $\chi:= \psi \circ \om \in
L^{2}(\nu)$ satisfies the $\si$-ASIP.
\end{rem}
We shall use Philipp-Stout's theorem (\cite{PS}, Section
7) to prove the ASIP for $\chi
:= \psi \circ \om$ on the Bernoulli space
$(\Sigma,s,\nu)$.  The version below comes from the original one
by using the $s$-invariance of $\nu$ and the independence of the
random process $(\xi_n)_{n \geq 0}$ defined by $\xi_n(\tal) = \al_n$.
\begin{theorem}{(Philipp-Stout)} Let $\chi$ be a $\nu$-centered observable on $\Sigma$ satisfying:
\begin{enumerate}
\item $\chi \in L^{2+\delta} (\nu)$ for some $\delta > 0$,
\item $\norm{  \chi - \E(\chi \vert \CC_n)  }_{2+\delta} \leq c \, \beta ^n$ for some $c > 0$ and $\beta < 1$.
\end{enumerate}
Then the sequence ${1 \over \sqrt n} \norm{S_n(\chi)}_2$ has a limit
$\si$. If $\si
> 0$, then $\chi$ satisfies the $\si$-ASIP.
\end{theorem}

\noindent Let us compare that result with Ibragimov's theorem (see \cite{I},
Theorem 2.1), which only requires moments of order $2$ and a summability condition:
\begin{theorem}{(Ibragimov)} Let $\chi$ be a $\nu$-centered observable on $\Sigma$
  satisfying:
\[ \sum_{n \geq 0}
  \norm{  \chi - \E(\chi \vert \CC_n)  }_2 < \infty . \]
Then the sequence ${1 \over
    \sqrt n} \norm{S_n(\chi)}_2$ has a limit
$\si$. If $\si > 0$, then $\chi$ satisfies the $\si$-CLT.
\end{theorem}

\section{A quantified version of Briend-Duval theorem}\label{BD}

This section is devoted to the proof of theorem \ref{good} (see subsection \ref{cons}). That
result will be crucial
to establish theorem A.

\subsection{Briend-Duval theorem}\label{cons}

We recall that $\VV_l = \cup_{i = 1}^ l f^i(\CC)$, $\VV = \cup_{i = 0}^\infty
f^i(\CC)$ and that $d_t = d^k$ is the
topological degree of $f$ (see subsection \ref{consa}). We set $\tau_*:= 2 \, \degr \VV_1 / (1-1/d)$.
\begin{theorem}{(Briend-Duval \cite{BrDu2})}
Let $\eta > 0$ and $l \geq 1$ be such
that $\tau_* / d^l < \eta$. Let $L$ be a complex line in $\Pj^k$ not contained in $\VV$, and
$\Delta \Subset \tilde \Delta$ be topological discs in $L \setminus \VV_l$. Then, for any $n \geq
l$, there exist $(1-\eta)d_t^n$ inverse branches $g_n$ on $\Delta$ satisfying:
\[  \Diam g_n (\Delta)  \leq { \tilde c \, d^{-n/2}   \over \eta^{1/2} \  \Mod (\tilde \Delta \setminus \Delta )^{1/2} }, \]
where $\tilde c$ is a universal constant, and $\Mod (\tilde \Delta \setminus \Delta )$ is the modulus of the annulus $\tilde \Delta \setminus \Delta$.
\end{theorem}
Let us recall the definition of the modulus (see Ahlfors book
\cite{A}, chapters
1 and 2). Let $\Lambda$ denote the family of
curves joining the boundary components of $A:= \tilde \Delta
\setminus \Delta$. For any conformal metric $\rho$ on $A$, we
respectively denote
 by $\area_\rho$ and by $l_\rho$ the area and the length
with respect to $\rho$. We denote by $\conf (A)$ the set of conformal
metrics giving finite area to $A$. The modulus of the annulus $A$ is then defined by:
\[ \Mod ( A ):= \sup_{\rho \in \conf (A)} \,   { l_\rho(\Lambda)^2 \over \area_\rho (A) } , \]
where $l_\rho(\Lambda):= \inf_{\lambda \in \Lambda} l_\rho
(\lambda)$.

\subsection{Statement of the quantified version}\label{cons}

We begin with some notations. Let $0 < \theta <1$ and $\theta_n:= [  \theta  n + { \log
    \tau_* \over  \log d } ] +1$. We introduce this integer in view
    of applying Briend-Duval theorem  with $\eta = d^{-\theta n}$
    and $l = \theta_n$ (indeed, $\tau_* / d^{\theta_n} < d^{-\theta
    n}$). Since the degree of $\VV_{\theta_n} = \cup_{i = 1}^{\theta_n} f^i(\CC)$ is at most $\tau_{\theta_n} = (d + \ldots + d^{\theta_n}) \deg \CC$, we have $\tau_{\theta_n} < d^{\theta_n}$ up to a multiplicative constant. \\

We let $0 < \theta < \theta' < 1$ and consider $n_0 \geq 1$ satisfying:
\begin{equation}\label{nzero}
\forall n \geq n_0  \ \ \ ,  \ \ \ \theta_n < \theta'n \ \ \textrm{ and }  \ \ \ \tau_{\theta_n} < {d}^{\, \theta' n}.
\end{equation}
Let us recall that $\VV_{\theta_n}[\delta]$ is the
$\delta$-neighbourhood of $\VV_{\theta_n}$ in $\Pj^k$. We fix
$\theta'/2 < \zeta < 1$ and define $\DD:= \limsup_{n \geq n_0}
\VV_{\theta_n}[d^{-\zeta n}]$.
\begin{prop}\label{vold}
The set $\DD$ satisfies $\Vol(\DD) = 0$.
\end{prop}
The proof is
postponed to subsection \ref{volo}. We now state the quantified version. The constant $\tilde c$ has been introduced in the statement of Briend-Duval's theorem, and we denote by $L$ the complex
line containing $z$ and $w$.
\begin{thm}\label{good}
There exists $\epsi > 0$ such that for every distinct points $(z,w) \notin \DD \cup \VV$, there exist an injective smooth path $\ga :
[0,1] \to L \setminus \VV$ joining $z$ and $w$, a decreasing family of
topological discs $(\Delta_n)_n \subset L$ and an integer $n_{z,w}$ such that for any $n \geq n_{z,w}$:
\begin{enumerate}
\item $\ga[0,1] \subset \Delta_n \subset L \setminus \VV_{\theta_n}$,
\item there exist $(1-d^{-\theta n }) d_t ^n$ inverse branches of $f^n$ on $\Delta_n$,
\item these branches satisfy $\Diam g_n (\Delta_n)  \leq \tilde c \,  d^{- \epsi n }$.
\end{enumerate}
\end{thm}
We note that $\theta,\epsi$ and $\tilde c$ do not depend on $(z,w) \in \Pj^k \setminus (\DD \cup \VV)$.

\subsection{Construction of good paths in the complex line $L \subset \Pj^k$} \label{ppaa}

Let $(z,w)$ be distinct points in $\Pj^k \setminus (\DD \cup \VV
)$. We identify the complex line $L$ containing $z$ and $w$ with
the $2$-dimensional sphere. We recall that the Fubini-Study metric induces on $L$ the standard
spherical metric $\sph$ with area $1$. We assume with no loss of generality
that $z$ and $w$ are the North and South pole of $L$. Let $E$ be the
equator of $L$. For any $y \in E$, we denote by $M_y$ the meridian
containing $y$, and by $M_y \{ \delta \}$ the $\delta$-neighbourhood of
$M_y$ in $L$ for the spherical metric. The constants $0 < \theta < \theta'
< 2\zeta$ have been defined in subsection \ref{cons}. Now we let
$0 < \zeta < \zeta' < \zeta'' < 1$ satisfying:
\begin{equation}\label{sop}
\theta' < \zeta'' - \zeta'   \ \textrm{ and } \  \theta + \zeta'' < 1.
\end{equation}
We may take for $(\theta,\theta',\zeta, \zeta' ,\zeta'')$ suitable
multiples of a small $\theta > 0$. The second inequality
of (\ref{sop}) will be used in next subsection. The integer $n_0$
has been defined in subsection \ref{cons}.

\begin{prop}\label{touch}
Let $(z,w)$ be distinct points in $\Pj^k \setminus (\DD \cup \VV )$. With the above notations, there exists a subset $F \subset E$ of full Lebesgue measure satisfying
the following properties. For any $y \in F$, there exists $n_{z,w}(y) \geq n_0$ such that:
\begin{enumerate}
\item the meridian $M_y$ does not intersect $\VV$,
\item the neighbourhood $M_y \{ d^{- \zeta'' n} \}$ does not intersect
  $\VV_{\theta_n}$ for any $n \geq n_{z,w}(y)$.
\end{enumerate}
\end{prop}
Let us now prove proposition \ref{touch}. We start with some
notations. Let $H^+$ and $H^-$ be the (open) North and South
hemispheres of $L$, these sets induce a partition $L = H^+ \sqcup E
\sqcup  H^-$. We denote by $\Leb$ the Lebesgue measure on $E$ and by
$p_1$ (resp. $p_2$) the spherical projection from $z$ (resp. $w$) to
$E$. For any $y \in E$ and $\delta > 0$, let $\II(y,\delta)$ be the
interval in $E$ centered at $y$ with length $2\delta$. We also
denote by $D(c,\delta) \subset L$ the disc with center $c$ and
radius $\delta$. We
  define $p_\kappa(c):= p_1(c)$ if $c \in H^+ \cup E$ and
  $p_\kappa(c):= p_2(c)$ if $c \in H^-$. The same convention holds
  for the projection of $D(c,\delta)$ to $E$: we use $p_1$ or $p_2$ depending on $c \in H^+ \cup E$ or $c \in H^-$.\\

Let $\{ c_i \, , \, 1 \leq i \leq
l_{\theta_n} \}:= \VV_{\theta_n} \cap L$, where $l_{\theta_n} \leq
\degr (\VV_{\theta_n}) \leq \tau_{\theta_n}$. Since the Fubini-Study metric induces $\sph$ on $L$, the set $\LL_{\theta_n}:=
  \cup_{i=1}^{l_{\theta_n}} D(c_i, d^{-\zeta n})$ is a subset of $\VV_{\theta_n}[ d^{-\zeta n} ]$. We recall that $\DD = \limsup_{n \geq n_0}
\VV_{\theta_n}[d^{-\zeta n}]$ and that $(z,w)
\notin \DD$. Thus there exists $n_1 \geq n_0$ depending on $(z,w)$ such that:
\begin{equation}\label{fl}
 \forall n \geq n_1 \ , \ (z , w)  \notin \VV_{\theta_n}[d^{-\zeta n}].
\end{equation}
In particular $(z,w) \notin
  \LL_{\theta_n}$. Since $\zeta < \zeta' < \zeta''$, we may increase $n_1$ so that $d^{- \zeta' n} + d^{- \zeta'' n} < d^{-\zeta n}$ for any $n \geq n_1$. We have therefore, for $\rho = z$ or $w$:
\begin{equation*}\label{evit}
  \forall \, 1 \leq i \leq l_{\theta_n} \ , \  \forall n \geq n_1  \ , \   D(\rho ,d^{-\zeta' n}) \cap
  D(c_i,d^{-\zeta'' n}) = \emptyset.
\end{equation*}
This implies, with $e_i:= p_\kappa(c_i) \in E$ and $c$ a positive
constant:
\begin{equation}\label{kkkk}
 \forall \, 1 \leq i \leq l_{\theta_n} \ , \ p_\kappa  \left( D (c_i,d^{-\zeta''n} ) \right) \subset \II_i:= \II(e_i \, , \,  c \,  d^{-\zeta'' n} \cdot d^{\zeta' n}).
\end{equation}
Hence $\II(\theta_n):= \cup_{i = 1}^{l_{\theta_n}} \II_i$ satisfies $\Lebe \II(\theta_n) \leq \tau_{\theta_n}  \cdot c \, d^{-(\zeta''- \zeta') n} \leq c \, d^{(\theta'- (\zeta'' - \zeta'))n}$. Since $\sum_n \Lebe \II(\theta_n) < \infty$ (see (\ref{sop})), the Borel-Cantelli lemma yields, for every $y$ in a full Lebesgue measure subset $F' \subset E$, an integer $n_{z,w}(y) \geq n_1$ satisfying:
\begin{equation}\label{flf}
y \notin  \bigcup_{n \geq n_{z,w}(y)} \II(\theta_n) .
\end{equation}
Let us prove the point 2 of proposition \ref{touch} (the point 1 will be proved
below, $F$ is a
subset of $F'$). Let $y \in F'$ and $\II:= \II(y,d^{ - (\zeta'' -
  \zeta') n})$. Since the intervals $\II_i$ defining $\II(\theta_n)$
are centered at $e_i = p_\kappa(c_i)$, the set $p_1^{-1}(\II)$ does
not intersect any point $c_i \in H^+ \cup E$. The same property holds for
$p_2^{-1} (\II)$ with the $c_i \in H^-$. This implies that $M_y \{
d^{-\zeta'' n} \}$ does not intersect $\VV_{\theta_n} \cap L$ for any
$n \geq n_{z,w}(y)$, and yields the point 2.

For the point 1, it suffices to verify that $p_\kappa(\VV \cap L)$ has zero Lebesgue
measure. Let $\WW:= \VV \cap L$. Since  $(z,w) \in L$ and $(z,w) \notin  \VV = \cup_{i =
  0}^\infty f^i(\CC)$, the complex line $L$ is not an algebraic
subset of the hypersurface $f^i(\CC)$ for any $i \geq 0$. In particular, $\WW_i:=  f^i(\CC) \cap L$ is finite for every $i \geq 0$. Hence  $\WW = \cup_{i \geq
  0} \WW_i$ satisfies $\Leb(p_\kappa(\WW)) = 0$. We finally set $F:= F' \setminus
p_\kappa(\WW)$, that completes the proof of proposition \ref{touch}.

\subsection{Proof of theorem \ref{good}} \label{bpo}

We set $\epsi:= {1 \over 2} (1 - (\theta + \zeta'')) > 0$ (see (\ref{sop})). Let $(z,w)$ be distinct
points in $\Pj^k \setminus (\DD \cup \VV )$ and consider some $y \in F$ provided by proposition \ref{touch}: the meridian $M_y$ does not intersect $\VV$ and its
neighbourhood $M_y \{ d^{- \zeta'' n} \}$ in $L$ does not intersect $\VV_{\theta_n}$
for every $n \geq n_{z,w}(y)$.  \\

We set $n_{z,w}:= n_{z,w}(y)$ and denote $M:= M_y$ for sake of simplicity. Let $\ga : [0,1] \to L$ be the natural parametrization of $M$. We define $\Delta_n:= M \{ d^{-\zeta'' n} / 2 \}$ and $\tilde \Delta_n:= M \{ d^{- \zeta'' n} \}$. Let us apply Briend-Duval's theorem with $\eta = d^{-\theta n}$, $l =
\theta_n$ and $\Delta_n  \Subset \ \tilde \Delta_n \subset \ L
\setminus \VV_{\theta_n}$. Since $n > \theta' n \geq \theta_n = l$ and
$\tau_* / d^{\theta_n} < d^{-\theta n}$ (see (\ref{nzero})), there exist $(1-d^{-\theta  n }) d_t^n$ inverse branches on the disc $\Delta_n$ satisfying:
\begin{equation}\label{kopa}
  \Diam g_n (\Delta_n)  \leq  \tilde c \,  d^{-n/2}  \, \left(  d^{-\theta n} \, \Mod \left[ \tilde \Delta_n \setminus \Delta_n \right] \right)  ^{-1/2} .
\end{equation}
It remains to bound the modulus of $A_n:= \tilde
\Delta_n \setminus \Delta_n$. Let $\Lambda_n$ be the set of curves
joining the boundary components of $A_n$. We denote by $\area_\sph$ and by $l_\sph$ the area and the length in $L$ with respect to the spherical metric $\sph$. The following estimates hold up to multiplicative constants. We have $l_\sph (\lambda) \geq d^{- \zeta'' n}$ for any $\lambda \in \Lambda_n$, hence $l_\sph (\Lambda_n) = \inf_{\lambda \in \Lambda_n} l_\sph
(\lambda) \geq  d^{- \zeta'' n}$. The inequalities $\area_\sph (A_n) \leq \area_\sph (\tilde \Delta_n) \leq d^{- \zeta'' n}$ then imply:
\begin{equation}\label{kopb}
\Mod ( A_n ) =   \sup_{\rho \in \conf A_n}  { l_\rho(\Lambda)^2 \over \area_\rho (A_n) } \geq  { l_\sph (\Lambda_n)^2 \over \area_\sph (A_n) } \geq { d^{- 2 \zeta'' n} \over d^{- \zeta'' n} } = d^{- \zeta'' n} .
\end{equation}
From (\ref{kopa}), (\ref{kopb}) and $\epsi = {1 \over 2 } (1 - (\theta + \zeta''))$,
we deduce that $\Diam g_n (\Delta_n) \leq \tilde c \, d^{-\epsi n}$. That completes the proof of theorem \ref{good}.

\subsection{Volume of neighbourhoods}\label{volo}

This subsection is devoted to the proof of proposition
\ref{vold}: we want to show
$\Vol(\DD) = 0$, where $\DD = \bigcap_{n \geq n_0} \bigcup_{p \geq
  n} \VV_{\theta_p}[d^{-\zeta p}]$. We recall that $\VV_{\theta_p}[d^{-\zeta
    p}]$ is the $d^{-\zeta p}$-neighbourhood of $\cup_{i =
  1}^{\theta_p} f^i(\CC)$ and that $\zeta > \theta'/2$. The proof is
based on the following lemma (see \cite{DS4}, lemma 2.3.8).

\begin{lem}\label{vois}
Let $X \subset \Pj^k$ be an algebraic subvariety of dimension $m$ and
degree $q$. Then $\Vol  X[\delta]  \leq q \, \delta^{2(k-m)}$ for any
$\delta > 0$, up to a multiplicative constant independent of $X$.
\end{lem}

We deduce $\Vol(\DD) = 0$ as follows. We set $p \geq n \geq n_0$ and
apply lemma \ref{vois} with $X = \VV_{\theta_p}$ and $\delta =
d^{-\zeta p}$ (here $k-m = 1$ and $q = \degr \VV_{\theta_p} \leq
\tau_{\theta_p}$). We obtain with $\tau_{\theta_p} \leq d^{\theta' p}$
(see (\ref{nzero})): $\Vol  \VV_{\theta_p}[d^{-\zeta p}]
 \leq \tau_{\theta_p} (d^{-\zeta p})^2  \leq d^{- (2\zeta -
   \theta')p}$. Hence:
\[ \forall n \geq n_0 \ , \ \Vol (\DD) \leq  \Vol  \bigcup_{p \geq n}  \VV_{\theta_p}[d^{-\zeta
    p}] \leq c_{2\zeta - \theta'} \, d ^{- ( 2\zeta - \theta' ) n} .  \]
This yields $\Vol (\DD) = 0$ when $n$ tends to infinity.

\proo \textsc{of lemma  \ref{vois}}: The argument is based on Lelong's inequality. Let $\EE$ be a \emph{maximal} $\delta$-separated
set in $X$ for the ambient metric: this means that $d(a,b) \geq \delta$ for any pair of distinct elements
of $\EE$, and that for any $x \in X$ there exists $a \in \EE$ satisfying
$d(a,x) < \delta$.  Since $X[\delta] \subset \cup_{a \in \EE}
B_a(2\delta)$, we get up to a multiplicative constant:
\begin{equation}\label{seff}
\Vol X[\delta]   \leq  (2\delta)^{2k} \, \cdb \EE .
\end{equation}
We now give an upper bound for $\cdb \EE$. Observe that $\Vol X$ is equal to the degree of $X$, and that the balls $(B_a (\delta / 2 ) )_{a \in \EE}$ are mutually disjoint. Thus:
\begin{equation*}\label{sefa}
  q = \Vol X \geq \sum_{a \in \EE} \Vol  \left( X \cap B_a (\delta / 2 ) \right).
\end{equation*}
Now Lelong's inequality asserts that $\Vol  (X \cap B_a (\delta /
2 )) \geq \delta^{2m}$ for any $a \in \EE$, up to a multiplicative
constant. Hence $\cdb \EE \leq q \, \delta^{-2m}$, as desired. \finsec

\section{Proof of theorem A}\label{GCT}

We set $\SS:= \VV \cup \DD \cup f(\DD) \cup
\per(f)$, where $\DD$ is defined in subsection
\ref{cons}. We have $\Vol(\SS)=0$ since $\Vol(\DD) = 0$. Let us recall the statement of theorem A.

\begin{theorem}{A:}
Let $z \in \Pj^k \setminus \SS$. There exist an $s$-invariant set $\Sigma'
\subset \Sigma$ of full $\nu$-measure and an $f$-invariant set $\JJ'
\subset \JJ$ of full $\mu$-measure satisfying the following properties. For any $\tal \in
\Sigma'$, the point $\om(\tal):= \lim_{n \to \infty} z_n(\tal) \in \JJ'$ is
well defined. We have $\om_* \nu  = \mu$ and the following diagram commutes:
\[ \xymatrix{
           \Sigma'  \ar[rr]^{s}  \ar[d]_\om &    & \Sigma'  \ar[d]^\om \\
           \JJ'  \ar[rr]^{f}       & &  \JJ'        }
\]
Moreover there exist $\theta , \epsi > 0$, $n_z \geq 1$ and $\tilde n: \Sigma' \to \N$ larger than $n_z$ such that:
\begin{enumerate}
\item $d(z_n(\tal),\om(\tal) ) \leq \tilde c_\epsi \, d^{-\epsi n}$  for every $\tal \in \Sigma'$ and $n \geq \tilde n (\tal)$,
\item $\nu ( \{ \tilde n  \leq q \} )  \geq 1 - c_\theta \, d^{-\theta q}$ for every $q \geq n_z$.
\end{enumerate}
\end{theorem}

We shall use theorem
\ref{good} and the method of coding trees introduced in \cite{PUZ} for $(\Pj^1,f,\mu)$. We recall
that $\AA = \{ 1, \ldots, d_t \}$. Let $z \notin \SS$ and $\{ w_\al \, , \, \al \in \AA \}:= f^{-1}(z)$. By the very definition of $\SS$, the cardinal of $f^{-1}(z)$ is equal
to $d_t$ and $w_\al \neq z$, $w_\al
\notin \VV \cup \DD$ for every $\al \in \AA$. We denote by
$L_\al$ the projective line in $\Pj^k$ containing $(z,w_\al)$ and
apply theorem \ref{good}: let $\ga_\al$ be an injective smooth
path joining $(z, w_\al)$ and $(\Delta_n(\al))_n \subset L_\al$ be a decreasing
sequence of discs containing $\ga_\al$ provided by that
theorem. We set $n_z:= \max \{ n_{z,w_\al} \, , \, \al \in \AA \}$. \\

Let us fix $\tilde \al = (\al_n)_{n \geq 0} \in \Sigma$. We define inductively injective smooth
paths $\ga_n(\tal) : [0,1] \to \Pj^k \setminus \VV$ and points
$z_n(\tal) \in \Pj^k \setminus \VV$. We first set $\ga_0 (\tal):=
\ga_{\al_0}$. This path joins $z =
\ga_0 (\tal)(0)$ and $w_{\al_0} = \ga_0 (\tal)(1) =: z_0
(\tal)$. Assume that the paths $\ga_j (\tal)$ and the points $z_j (\tal)$ have
been defined for $0 \leq j \leq n-1$. We let $\ga_n(\tilde \al)$ to be
the lift of $\ga_{\al_n}$ by $f^n$ with starting point $\ga_n(\tilde
\al)(0) = z_{n-1} (\tal)$. This path is well defined since $\ga_{\al_n}$ does
not intersect $\VV$. We finally let $z_n(\tal):=
\ga_n (\tal)(1)$. \\

We note that $z_{n-1}(\tal)$ and
$z_n(\tal)$ are the endpoints of $\ga_n(\tal)$ and that
$z_n(\Sigma) = f^{-(n+1)}(z)$ has cardinal
$d_t^{n+1}$. The reader will easily check
the relation $f \circ z_n(\tal) = z_{n-1} \circ s (\tal)$. Observe
also that $\ga_n(\tal)$ and $z_n(\tal)$ depend only on
$\pi_n(\tal) = (\al_0,\ldots,\al_n)$. The following lemma is a
consequence of theorem \ref{good} and the fact that $\ga_\al[0,1] \subset
\Delta_n(\al)$.

\begin{lem}\label{decr}
For every $\al \in \AA$ and $n \geq
n_z$, there exist at least $(1-d^{-\theta n})d_t^n$ elements $(\al_0,\ldots,\al_{n-1}) \in \AA^n$ such that $\Diam  \ga_n (\al_0,\ldots,\al_{n-1},\al)  \leq \tilde c \,  d^{- \epsi n }$.
\end{lem}
Let $\Omega_n:= \{ \tilde \al \in \Sigma \, , \, \Diam \ga_n(\tal) > \tilde c \,
d^{-\epsi n} \}$ and $\BB_n$ be the collection of $(n+1)$-cylinders $\{ C_n(\tal)  \ , \ \tal \in \Omega_n \}$. We have $\Omega_n = \BB_n^*$. Let us also  define:
\[  \Omega(n) \ := \ \bigcup_{p \geq n} \Omega_p \  =   \  \bigcup_{p \geq n} \BB_p^* \, .  \]

\begin{lem}\label{cdi} For any $n \geq n_z$, we have:
\begin{enumerate}
\item $\cdb (\BB_n) \leq d_t^{n+1}  d^{-\theta n}$.
\item $\nu(\Omega_n) \leq d^{-\theta n}$, hence $\nu (\Omega(n))  \leq c_\theta \, d^{-\theta n}$.
\item if $\tal \notin \Omega(n)$, then $d( z_{m-1} (\tal), z_m(\tal)
  ) \leq \tilde c \, d^{-\epsi m}$ for any $m \geq n$.
\end{enumerate}
\end{lem}

\proof  We have $\BB_n = \{ C_n(\tal)  \ , \ \Diam \ga_n(\tal) > \tilde c \,
d^{-\epsi n}   \}$. For every $\al \in \AA$, we set $\BB_n(\al)
\subset \BB_n$ to be the collection of
$(n+1)$-cylinders whose last coordinate is equal to $\al$. The lemma
\ref{decr}  implies that $\cdb (\BB_n(\al))
\leq d_t^n  d^{-\theta n}$ and thus $\cdb (\BB_n) = \sum_{\al \in \AA} \cdb( \BB_n(\al)
) \leq d_t^{n+1} d^{-\theta n}$, which is the point 1. The point 2
follows:
\[ \nu(\Omega_n)
= \nu(\BB_n^*) = \cdb (\BB_n) / d_t^{n+1} \leq d^{-\theta
  n}. \]
For the point 3, observe that $d(  z_{m-1} (\tal) , z_m(\tal) ) \leq
  \Diam \ga_m(\tal)$. If $\tal \notin \Omega(n)$, then $\tal \notin
  \Omega_m$ for any $m \geq n$, hence $\Diam \ga_m(\tal) \leq
\tilde c \, d^{-\epsi m}$. \fin

Let $\Omega:= \bigcap_{n \geq n_z} \Omega(n) =  \limsup_{n
  \geq n_z} \Omega_n$. The set $\Sigma'':= \Sigma \setminus
\Omega$ has full $\nu$-measure since $\nu(\Omega) \leq \nu(\Omega(n)) \leq c_\theta d^{-\theta n}$ for any $n
\geq n_z$. For every $\tal \in \Sigma''$, we define $\tilde n (\tal)$ to be the least integer $n \geq n_z$ satisfying $\tal \notin
\Omega(n)$. Let $\Theta_q:= \{ \tilde n  \leq q \}$.

\begin{lem}\label{pushoo} $ $
\begin{enumerate}
\item $\om(\tal)=\lim_{n
  \to \infty} z_n(\tal)$ is well defined for every $\tal \in \Sigma''$.

\item $d(
  z_n(\tal), \om(\tal) ) \leq \tilde c_\epsi \, d^{-\epsi n}$ for every $n
  \geq \tilde n (\tal)$.

\item $\om : \Sigma'' \to \Pj^k$ satisfies $\om_*\nu = \mu$.

\item $\nu ( \Theta_q )  \geq 1 - c_\theta \, d^{-\theta q}$ for any $q \geq n_z$.
\end{enumerate}
\end{lem}

\proof
The points 1, 2 and 4 come from lemma \ref{cdi}(3,2) and
the definition of $\tilde n (\tal)$. Now we prove the point 3. Let us consider the surjective map $z_n : \Sigma'' \to  f^{-(n+1)}(z)$. Since
$z_n(\tal)$ depends only on $\underline \al:=
(\al_0,\ldots,\al_n) \in \AA^{n+1}$, the measure  ${z_n}_* \nu$ is
equal to:
\[ {z_n}_* \nu = \sum_{\underline \al \in \AA^{n+1}} \nu \left( \Sigma''
\cap C_n( \underline \al ) \right) \, \delta_{z_n(\underline \al)} =
 {1 \over d_t^{n+1}} \, \sum_{f^{n+1}(y) = z} \delta_y = \mu_{n+1,z}. \]
Since $z \notin \SS$ and $\EE \subset \VV \subset
\SS$, the sequence of probability measures $(\mu_{n,z})_n$ converges to
$\mu$ (see subsection \ref{consa}). Hence it remains to prove ${z_n}_* \nu \to \om_* \nu$, meaning that
 $\int_{\Sigma''} \varphi \circ z_n \, d\nu \to \int_{\Sigma''}
 \varphi \circ \om \, d\nu$ for every test function $\varphi :
 \Pj^k \to \R$. But this follows from point 1 and Lebesgue convergence theorem.  \fin

It remains to define $\Sigma' , \JJ'$ and to verify
the relation $f \circ \om = \om \circ s$ on $\Sigma '$. The lemma \ref{pushoo}(3) implies that $\Sigma_*:= \om(\Sigma'')$
satisfies $\mu(\Sigma_*) = \nu(\om^{-1}\Sigma_*) \geq \nu
(\Sigma'') = 1$. We define $\JJ':=  \bigcap_{n
  \in \Z} f^n(\JJ \cap \Sigma_*) $  and $\Sigma':=   \bigcap_{n \in \Z}
s^n(\Sigma'' \cap \om^{-1} \JJ')$. These are invariant subsets of full
measure. We obtain $f \circ \om = \om \circ s$ on $\Sigma'$ by taking
limits in $f \circ z_n(\tal) = z_{n-1} \circ s (\tal)$. That completes the proof of theorem A.

\section{Proof of theorem B} \label{estiho}

Let us recall the statement.

\begin{theorem}{B:} Let $\psi \in \UU$ be a
  $\mu$-centered observable and $\om$ be a coding map provided by theorem A. Let $\chi:= \psi \circ
  \om$ and $1 \leq p < +\infty$.
\begin{enumerate}
\item there
  exist $\hat c_p, \lambda_p > 0$ such that $\norm { \chi   - \E ( \chi \vert
  \CC_n  )  }_p \leq \hat c_p \, e^{-n\lambda_p}$ for every $n \geq 0$.

\item $R_j(\chi):= \int_\Sigma \chi \cdot
\chi \circ s^j \, d\nu$ satisfies $\abs{R_j(\chi)} \leq 2 \norm {\chi}_2  \hat c_2 \,
  e^{-(j-1) \lambda_2}$ for every $j \geq 1$.
\end{enumerate}
\end{theorem}

\subsection{Proof of theorem B(1)}

We set $\chi_B:= \chi .  1_B$ for any $B
\subset \Sigma$ and use the following estimates provided by
theorem A. We recall that $\Theta_n = \{ \tilde n (\tal) \leq n \}$.
\begin{enumerate}
\item[($\star$)]  $d( z_n(\tal),
\om(\tal) ) \leq \tilde c_\epsi \, d^{-\epsi n}$ for every $\tal \in \Theta_n$,
\item[($\star \star$)] $\nu(\Theta_n) \geq 1 -  c_\theta  \, d^{- n
  \theta}$ for every $n \geq n_z$.
\end{enumerate}
We will need the following lemma, which is a direct consequence of ($\star$).
\begin{lem}\label{fcx}
Let $\tal \in \Theta_n$ and $\tbe \in C_n(\tal) \cap \Theta_n$. Then $d(\om(\tal),\om(\tbe)) \leq 2 \, \tilde c_\epsi \, d^{-\epsi n}$.
\end{lem}

\subsubsection{The H\"older case}\label{ccc2}

Let $\psi$ be an $h$-H\"older and $\mu$-centered observable on $\Pj^k$. We set $\chi:= \psi \circ
\om$. The theorem B(1) is a consequence of the following estimates, which
hold for every $n \geq n_z$.

\begin{lem}\label{esb1}
$\norm { \chi_{\Theta_n^c}   - \E ( \chi_{\Theta_n^c} \vert \CC_n  )
  }_p \leq 2 \norm{\chi}_\infty  \, \left( c_\theta   \, d^{- n  \theta} \right)^{1/p} $.
\end{lem}

\proof
The left hand side is less than $2 \norm {\chi_{\Theta_n^c} }_p$ by
Jensen inequality. Then the conclusion follows from $(\star \star)$. \finsec

\begin{lem}\label{esb2}
$\norm { \chi_{\Theta_n} - \E ( \chi_{\Theta_n} \vert \CC_n  )
  }_p \leq c \, d^{-n \tau}$ for some $c, \tau > 0$.
\end{lem}

\proof We denote $\varphi:= \chi_{\Theta_n} - \E ( \chi_{\Theta_n} \vert
\CC_n)$ and estimate  $\norm
   {\varphi_{\Theta_n^c}}_p$, $\norm{\varphi_{\Theta_n}}_p$. Since $\varphi_{\Theta_n^c} =  - \E ( \chi_{\Theta_n} \vert
\CC_n) \cdot 1_{\Theta_n^c}$, we have:
\[ \norm{\varphi_{\Theta_n^c}}_p \leq \norm { \E ( \chi_{\Theta_n} \vert
\CC_n)  }_{2p} \cdot \nu(\Theta_n^c) ^{1/2p} \leq \norm{\chi}_{2p} \cdot \left( c_\theta \, d^{- n  \theta} \right)^{1/2p}.\]
We now deal with $\norm{\varphi_{\Theta_n}}_p$. For every $\tal \in
\Theta_n$, let $\nu_\tal$ be the conditional
measure of $\nu$ on the cylinder $C_n(\tal)$. We have for
every $\tal \in \Theta_n$:
\begin{equation}\label{lki}
  \varphi_{\Theta_n}(\tal) = \int_{C_n(\tal) \cap \Theta_n} \left( \chi(\tal)
-\chi (\tbe) \right)  d\nu_\tal(\tbe)  + \chi(\tal) \cdot  \nu_\tal(C_n(\tal) \cap \Theta_n^c).
\end{equation}
We deduce from $\chi =  \psi \circ \om $, lemma \ref{fcx} and the
fact that  $\psi$ is $h$-H\"older:
\begin{equation*}
 \forall \tal \in \Theta_n \ , \ \abs{\varphi_{\Theta_n}(\tal)} \leq \left( 2 \, \tilde
 c_\epsi \, d^{-n\epsi} \right)^h  + \norm{\chi_{\Theta_n}}_\infty \cdot  \nu_\tal(C_n(\tal) \cap \Theta_n^c).
\end{equation*}
Hence we get for every $p \geq 1$ up to a multiplicative constant:
\begin{equation*}
 \forall \tal \in \Theta_n \ , \ \abs{\varphi_{\Theta_n}(\tal)}^p \leq d^{-nhp
 \epsi} + \norm{\chi_{\Theta_n}}^p_\infty  \cdot  \nu_\tal(C_n(\tal) \cap \Theta_n^c).
\end{equation*}
By integrating over $\Theta_n$ and using $(\star \star)$, we deduce:
\[ \norm{\varphi_{\Theta_n}}_p^p \leq d^{-nhp
 \epsi} +  \norm{\chi}^p_\infty  \cdot  c_\theta d^{- n
 \theta}. \]
That completes the proof of the lemma.  \finsec

\subsubsection{The general case $\psi \in \UU$} \label{estib}

Let $\psi : \Pj^k \to \R \cup \{ - \infty \}$ be a $\mu$-centered observable in
 $\UU$: the function $e^\psi$ is $h$-H\"older  and satisfies $\psi \geq \log
 d(\cdot,\NN_\psi)^\rho$ on $\Pj^k$ (see definition
 \ref{unb}). Observe in particular that $\psi$ is bounded from above. We recall that $\NN_\psi[r]$ is the
 $r$-neighbourhood of $\{ \psi = - \infty\}$ and that $\chi = \psi
 \circ \om$. We consider the following subsets of $\Sigma$:
\[ \SS_n:=
\Theta_n^c  \setminus \NN_n \ \ \ , \ \ \ \Gamma_n = \Theta_n \setminus \NN_n   \ \ \ , \ \ \      \NN_n:=  \om^{-1} \left( \NN_\psi[d^{-n (h \epsi / 2 \rho)}]
\right).  \]
We shall need the following observations. First, we have $\nu(\NN_n) = \mu ( \NN_\psi[d^{-n (h \epsi / 2 \rho)}] ) \leq d^{-n \ga (h \epsi / 2 \rho)}$ up to a multiplicative constant (see subsection \ref{consa}). We deduce from $(\star \star)$:
\begin{equation}\label{vatoo}
 \nu(\Ga_n^c) = \nu ( \Theta_n^c \cup \NN_n  ) \leq  c_\theta d^{- n
 \theta} + d^{-n \ga (h \epsi / 2 \rho)} \leq c \, d^{- n \eta }
\end{equation}
for some $c, \eta > 0$. Second, for every
  $\tal \in \NN_n^c = \SS_n \cup \Gamma_n$, we have $\chi(\tal) \geq \log d(\om(\tal) ,
  \NN_\psi)^\rho \geq \log d^{-\rho n (h \epsi / 2 \rho)}$, hence:
\begin{equation}\label{vato}
 \norm{\chi_{ \SS_n \cup \Gamma_n} }_\infty
  \leq n \, (h \epsi \log d ) / 2.
\end{equation}
The theorem B(1) is now a consequence of the following estimates.

\begin{lem}
$\norm {\chi_{\NN_n} - \E ( \chi_{\NN_n} \vert \CC_n ) }_p \leq \left( \kappa \, d^{- n (h \epsi / 2 \rho)\cdot (\ga/ 2) } \right)^{1/p}$.
\end{lem}

\proof
The left hand side is less than $2 \norm {
  \chi_{\NN_n} }_p$.  Proposition \ref{cuun} yields  $\norm {\chi_{\NN_n}}_p =  \norm {\psi \circ \om \cdot 1_{\NN_n}}_p  \leq \left( \kappa \,   d^{-n(h \epsi / 2
  \rho)\cdot(\ga /2)} \right)^{1/p}$ for every $n$ such that $d^{-n(h \epsi / 2
  \rho)} < 1/2$. \finsec

\begin{lem}
$\norm {\chi_{\SS_n} - \E ( \chi_{\SS_n} \vert \CC_n ) }_p \leq  n \, (h \epsi \log d) \cdot \left( c \, d^{- n  \eta} \right)^{1/p}$.
\end{lem}

\proof The left hand side is less than $2\norm {
  \chi_{\SS_n} }_p$. We
  conclude by using  (\ref{vatoo}) and (\ref{vato}) (observe that $\SS_n \subset \Gamma_n^c$).  \finsec

\begin{lem}
$\norm { \chi_{\Gamma_n}   - \E ( \chi_{\Gamma_n} \vert \CC_n  )
  }_p \leq c \, d^{-n\tau}$ for some $c, \tau > 0$.
\end{lem}

\proof We follow the proof of lemma \ref{esb2}: we set $\varphi:= \chi_{\Gamma_n} - \E ( \chi_{\Gamma_n} \vert
\CC_n)$ and estimate $\norm
   {\varphi_{\Gamma_n^c}}_p$, $\norm{\varphi_{\Gamma_n}}_p$. The line (\ref{vatoo}) yields:
\[ \norm{\varphi_{\Gamma_n^c}}_p \leq \norm { \E ( \chi_{\Gamma_n} \vert
\CC_n)  }_{2p} \cdot \nu(\Gamma_n^c) ^{1/2p} \leq \norm{\chi}_{2p}
\cdot \left( c \, d^{- n \eta} \right)^{1/2p}. \]
Now we deal with $\norm{\varphi_{\Gamma_n}}_p$. We can write as in (\ref{lki}):
\begin{equation}\label{popp}
 \forall \tal \in \Gamma_n \ , \ \varphi(\tal) = \int_{C_n(\tal) \cap \Gamma_n} \left( \chi(\tal)
-\chi (\tbe) \right)  d\nu_\tal(\tbe)  + \chi(\tal) \cdot  \nu_\tal(C_n(\tal) \cap \Gamma_n^c).
\end{equation}
Let $\tal \in \Gamma_n$ and $\tbe \in C_n(\tal) \cap \Gamma_n$. We deduce from $(\tal,\tbe) \notin \NN_n$ that $e^\psi \circ \om (\tal)$
 and $e^\psi
 \circ \om (\tbe)$ are larger than $d^{-n h \epsi  / 2 }$. This
 implies:
\begin{equation*}
\abs{ \chi(\tal) -\chi(\tbe) } = \abs{ \log e^\psi \circ \om (\tal)   -  \log e^\psi \circ \om
 (\tbe)}  \leq d^{n h \epsi / 2}  \abs{ e^\psi \circ \om (\tal)   -  e^\psi \circ \om (\tbe) }.
\end{equation*}
Using lemma \ref{fcx}
and the fact that  $e^\psi$ is
$h$-H\"older, the last term is less than $d^{n h \epsi / 2} \cdot
\left( 2 \, \tilde c_\epsi \, d^{-n \epsi}\right)^h$. Then we deduce from (\ref{popp}), up to a multiplicative constant:
\begin{equation*}
 \forall \tal \in \Gamma_n \ , \ \abs{\varphi(\tal)} \leq   d^{-n h \epsi / 2} + \norm{\chi_{\Gamma_n}}_\infty \cdot  \nu_\tal(C_n(\tal) \cap \Gamma_n^c).
\end{equation*}
Taking the $p$-th power, integrating over $\Gamma_n$ and using (\ref{vatoo}), (\ref{vato}), we obtain up to a multiplicative constant:
\[ \norm{\varphi_{\Gamma_n}}_p^p \leq d^{-nhp
 \epsi/2} +  \left( n ( h \epsi \log d ) / 2 \right)^p  \cdot  c \, d^{- n
 \eta} .\]
That completes the proof of the lemma. \finsec

\subsection{Proof of theorem B(2)}\label{corre}

Let $\psi \in \UU$ be a $\mu$-centered observable and $\chi
= \psi \circ \om$. Let $j \geq 1$ and $n \geq 0$ to be specified later. We
set $\chi_n:= \E(\chi \vert \CC_n)$ and write:
\[  \chi \cdot \chi \circ s^j = (\chi - \chi_n) \cdot \chi \circ s^j + \chi_n \cdot ( \chi \circ s^j - \chi_n \circ s^j ) + \chi_n \cdot \chi_n \circ s^j  .       \]
By using the $s$-invariance of $\nu$ and Jensen inequality $\norm{\chi_n}_2 \leq \norm{\chi}_2$, we deduce:
\begin{equation}\label{indp}
  \abs{R_j(\chi)}  = \left \vert  \int_{\Sigma} \chi \cdot \chi \circ s^j \, d\nu \right \vert  \leq 2  \norm{\chi}_2 \norm{\chi - \chi_n}_2 + \left \vert {\int_{\Sigma} \chi_n \cdot \chi_n \circ s^j \, d\nu } \right \vert .
\end{equation}
The variables  $\chi_n$ and $\chi_n
\circ s^j$ respectively depend
 on $(\xi_0,\ldots,\xi_n)$ and
$(\xi_j,\ldots,\xi_{n+j})$, where $\xi_n : \Sigma \to \AA$ is the
projection $\xi_n(\tal) = \al_n$. These are
independent variables when $n = j-1$, hence $\int_{\Sigma} \chi_n \cdot \chi_n \circ s^j \, d\nu = \int_{\Sigma}
\chi_n \, d\nu  \int_{\Sigma} \chi_n \circ s^j \, d\nu$ in that case. But this product is zero since $\chi$ is
$\nu$-centered. The conclusion then follows from (\ref{indp}) with $n = j-1$ and theorem B(1).

\section{Proof of theorem C} \label{CC}

Let us recall the statement.

\begin{theorem}{C:} For every $\mu$-centered observable $\psi \in \UU$, we have:
\begin{enumerate}
\item $\si:= \lim_{n \to \infty} {1 \over \sqrt n} \norm{S_n(\psi)}_2$ exists, and $\si^2 = \int_{\Pj^k} \psi^2 \, d\mu + 2\sum_{j \geq 1} \int_{\Pj^k}
 \psi \cdot \psi \circ f^j \, d\mu$.

\item If $\si = 0$, then $\psi = u - u \circ f$ $\mu$-a.e. for some $u \in L^2(\mu)$.

\item If $\si > 0$, then $\psi$ satisfies the $\si$-ASIP.
\end{enumerate}
\end{theorem}

\noindent The points 1 and 2 are consequences of classical lemma \ref{somme} below, whose condition $\sum_{j \geq 1} j \abs{R_j (\varphi)} < \infty$
is fulfilled by theorem B(2). The
point 3 follows from  proposition \ref{cuun}, theorem B(1) and Philipp-Stout's theorem (see subsection
\ref{consb}).

\begin{lem}\label{somme}
Let $(X,g,m)$ be a dynamical system and $\varphi \in L^2(m)$ be a
$m$-centered observable. We denote $S_n(\varphi) = \sum_{j = 0}^{n-1}
\varphi \circ g^j$ and $R_j(\varphi) = \int_X \varphi \cdot \varphi
\circ g^j \, dm$. Let $\si^2:=  R_0(\varphi) + 2\sum_{j \geq 1}
R_j(\varphi)$. If $\sum_{j \geq 1} j \abs{R_j (\varphi)} < \infty$, then
$\si^2$ is finite and we have:
\begin{enumerate}
\item $\norm{S_n(\varphi)}^2_2 = n \si^2 + O(1)$. In particular,
  $\lim_{n \to \infty} {1 \over n} \norm{S_n(\varphi)}^2_2 = \si^2$.
\item $\si^2 = 0$ if and only if $\varphi = u - u \circ g$ $m$-a.e. for some $u \in L^2(m)$.
\end{enumerate}
\end{lem}

\proof Let $S_n:= S_n(\varphi)$ and $R_j:= R_j (\varphi)$.  Since $m$ is $g$-invariant, we have $\norm{S_n}^2_2 = n R_0 + 2 \sum_{j=1}^{n-1} \left( n-j  \right)
R_j$. We deduce for every $n \geq 1$:
\begin{equation}\label{pml}
   \norm{S_n}^2_2 = n \left(  R_0 + 2 \sum_{j = 1}^\infty R_j
\right) +  (- 2) \left( \sum_{j = 1}^{n-1} j R_j +   \sum_{j = n}^\infty n R_j\right) = n \si^2 + A_n ,
\end{equation}
where $\abs{A_n} \leq  2 \sum_{j \geq 1} j \abs{R_j}$. That proves the point 1. Let us show the point 2. Suppose $\si^2 =
0$. In view of (\ref{pml}), the function $u_p:= {1 \over p} \sum_{n=1}^p S_n$ satisfies $\norm{u_p}_2 \leq ( 2 \sum_{j \geq 1} j \abs{R_j} )^{1/2}$ for every $p \geq 1$. Let $u:= \lim_{j \to \infty} u_{p_j}$ be a weak cluster point in $L^2(m)$ and observe that:
\[ \forall j \geq 1 \ , \ u_{p_j} - u_{p_j} \circ g = {1 \over p_j}  \sum_{n=0}^{p_j - 1} \left( \varphi - \varphi \circ g^n \right) = \varphi -  {1 \over p_j}  S_{p_j} .\]
We deduce $\varphi = u - u \circ g$ $m$-a.e. by taking limits in
$L^2(m)$ : $\lim_{j \to \infty}
u_{p_j} \circ g = u \circ g$ since $m$ is $g$-invariant, and $\lim_{j \to
  \infty}{1 \over p_j} S_{p_j} = \int_X \varphi \, dm = 0$ by Von
Neumann theorem. The reverse implication of the point 2 comes from
$\si^2 = \lim_{n \to \infty} {1 \over n} \norm{S_n(\varphi)}^2_2 =
\lim_{n \to \infty} {1 \over n} \norm{u - u \circ g ^n}^2_2 = 0$. \finsec

\section{Proof of theorem D} \label{appll}

We recall that $J:= \log \Jac f - \int_{\Pj^k} \log \Jac f \, d\mu$,
this is an unbounded $\mu$-centered observable in $\UU$. We set
$\si_J:= \lim_n {1 \over \sqrt n} \norm{S_n(J)}_2$, which is well
defined by theorem C. We denote by  $\chi_1 \leq \ldots \leq \chi_k$
the Lyapunov exponent of $\mu$, they are larger than or equal to
$\log d^{1/2}$.

\begin{theorem}{D:}
If the Lyapunov exponents of $\mu$ are minimal equal to $\log
d^{1/2}$, then $\si_J = 0$ and $\mu$ is absolutely continuous with
respect to Lebesgue measure.
\end{theorem}

The first part $\si_J = 0$ will be proved in subsection
\ref{zeze}. The second part is a consequence of theorem \ref{laa}
below (that theorem will be proved in subsection \ref{mpo} by using $\si_J = 0$). In the sequel, the maps $f^n$
and $d_x f^n$ are implicitely written in some fixed charts of $\Pj^k$.

\begin{thm}\label{laa}
Assume that the Lyapunov exponents are minimal. Then for $\mu$
almost every $x \in \Pj^k$, there exists $\rho(x) > 0$ and a
subsequence $(n_j(x))_{j \geq 1}$ such that $f^{n_j} \circ \left(x +
d^{-n_j/2} \cdot \Id_{\C^k} \right) : B(\rho(x)) \to \Pj^k$ is
injective.
\end{thm}

\proo \textsc{of the second part of theorem D (abolute continuity): } We use the notations of theorem \ref{laa}. Let $x \in
\Pj^k$ be a $\mu$-generic point and set $n_j:= n_j(x)$. Since
$f^{n_j}$ is injective on the ball $B_j:= B_x ( \rho(x)d^{-n_j/2} )
$ and $\mu$ has constant jacobian $d^k$ (see subsection
\ref{consa}), we obtain $\mu (B_j) =  \mu (f^{n_j} (B_j)) d^{-k n_j}$. Observe also that $\Leb (B_j) = \rho(x)^{2k} \left(
d^{-n_j / 2}  \right)^{2k} =  \rho(x)^{2k} d^{-k n_j}$ up to a
multiplicative constant. We obtain therefore for $\mu \ale \, x \in
\Pj^k$:
\[  \liminf_{r \to 0}  { \mu (B_x(r)) \over \Leb (B_x(r)) }  \leq \liminf_{j \to \infty} { \mu (B_j) \over \Leb (B_j) } =  \liminf_{j \to \infty} {\mu ( f^{n_j} (B_j) ) \over \rho(x)^{2k} } \leq {1 \over  \rho(x)^{2k} } < \infty. \]
That proves the absolute continuity of $\mu$ (see \cite{Mat}, theorem 2.12). \finsec

\subsection{Preliminaries} \label{abvc}

Observe that $J = \log \Jac f - \log d^k$ when the Lyapunov
exponents are equal to $\log d^{1/2}$. Since the jacobian is a
multiplicative function, we have in that case:
\begin{equation}\label{jjaa}
 S_n(J) = \sum_{i=0}^{n-1} J \circ f^i(x) = \log  \Jac f^n - \log d^{kn}.
\end{equation}
The \emph{singular values}  $\delta_1 \leq
\ldots \leq \delta_k$ of the linear map $A:= d_x f^n$ are defined as the
eigenvalues of $\sqrt{A A^*}$. In particular, there exist unitary
matrices $(U,V)$ such that $d_x f^n = U \,
\Diag(\delta_1,\ldots,\delta_k) \, V$. We have therefore:
\begin{equation}\label{delja}
\delta_1 = \norm {(d_x f^n)^{-1} }^{-1} \ \  \textrm{ and }  \ \  \prod_{i=1}^k {\delta_i}^2 = \Jac f^n(x) \geq \delta_1^{2k}.
\end{equation}
For any $\rho,\tau > 0$ and $n \geq 1$, we define:
\[  \BB_n(\rho):= \left \{  x \in \Pj^k \, , \, f^n \circ (x + d_x f^n)^{-1}  : B(\rho) \to \Pj^k \textrm{ is an injective map} \right\} ,  \]
\[  \RR_n(\tau):= \left\{ x \in \Pj^k \, , \, \norm {(d_x f^n)^{-1} }^{-1}  \geq d^{n/2}/\tau  \right \} .\]
The following estimates were proved by Berteloot-Dupont \cite{BeDu}.
They hold for every system $(\Pj^k, f, \mu)$ whose Lyapunov
exponents satisfy $\chi_k < 2 \chi_1$.
\begin{thm}\label{bdu}  There exists $\al : \, ]0,1] \to \R^*_+$
      satisfying $\lim_{\rho \to 0} \al(\rho) = 1$ and for $n \geq
      1$:
\begin{enumerate}
\item $\mu(\BB_n(\rho)) \geq \al(\rho)$,
\item  $\mu(\BB_n(\rho) \cap \RR_n(\tau)^c ) \leq (\rho \, \tau)^{-2}$.
  \end{enumerate}
\end{thm}
That result implies the following lemma.

\begin{lem}\label{eoui}
Let $\rho \in ]0,1]$. There exists $\HH \subset \Pj^k$ satisfying $\mu(\HH)=1$ and:
\[ \forall x \in \HH \ , \ \exists n(x) \geq 1 \ , \ \forall n \geq
n(x) \ , \ x \notin \BB_n(\rho)  \textrm{ or }  \Jac f^{n}(x) \geq
d^{kn} / n^{2k} . \]
\end{lem}

\proof We apply proposition \ref{bdu}(2) with $\tau = n$ to get $\mu (\BB_n(\rho) \cap \RR_n(n)^c) \leq  (\rho \, n)^{-2}$. Since $\sum_{n
  \geq 1} \mu (\BB_n(\rho) \cap \RR_n(n)^c) < \infty$, there exists by
  Borel-Cantelli lemma a subset $\HH$ of
  full $\mu$-measure satisfying:
\[ \forall x \in \HH \ , \ \exists n(x) \geq 1 \ , \ \forall n \geq
n(x) \ , \ x \notin \BB_n(\rho) \textrm{ or }  x \in  \RR_n(n) . \]
But $x \in  \RR_n(n)$ implies by (\ref{delja}):
$\Jac f^{n}(x) \geq \left( d^{n/2} / n \right)^{2k} = d^{kn} /
n^{2k}$ . \finsec

\subsection{Proof of the first part of theorem D ($\si_J = 0$)}\label{zeze}

Suppose that the exponents are minimal and that $\si_J  = \lim_n
{1 \over \sqrt n} \norm{S_n(J)}_2 > 0$. Then $J$ satisfies the CLT: if $V
:= {1 \over \sqrt{2\pi}} \int_{- \infty}^{-1}  e^{-u^2/2} \, du$, we
get $\mu \left( \GG_n:= \{  {S_n(J) \over \sqrt n } \leq - \si_J \} \right) \geq V / 2$ for $n$ larger than some $N$ (see subsection
\ref{consb}). \\

Let $\rho > 0$ be such
that $\mu(\BB_n(\rho)) > 1 - V / 4$ for every $n \geq 1$. If we set $\FF_n:= \BB_n(\rho) \cap \GG_n$, then $\FF
:= \limsup_{n \geq N} \FF_n$ satisfies $\mu(\FF)
\geq V/4$. Let $x \in \FF
\cap \HH$, where $\HH$ is provided by lemma \ref{eoui}. Let $(n_j(x))_{j \geq 1}$ be such that $x \in
\FF_{n_j}$ for every $j \geq 1$. The inclusion $\FF_{n_j}  \subset \GG_{n_j}$ yields
$S_{n_j}(J)(x) \leq - \si_J \sqrt{n_j}$ for every $j \geq 1$. Since $S_{n_j}(J) = \log \Jac f^{n_j} - \log d^{kn_j}$ (the exponents are indeed minimal, see (\ref{jjaa})), we deduce:
\begin{equation}\label{ma}
  \forall j \geq 1 \ , \ \Jac f^{n_j}(x) \leq d^{k n_j} e^{- \si_J \sqrt{n_j}}.
\end{equation}
But $\Jac f^{n_j}(x) \geq
d^{kn_j} / n_j^{2k}$ for every $n_j \geq n(x)$, following from $x \in \BB_{n_j}(\rho) \cap \HH$ and lemma \ref{eoui}. That contradicts (\ref{ma}) when $j$ tends to infinity.

\subsection{Proof of theorem \ref{laa}}\label{mpo}

We proved in subsection \ref{zeze} that $\si_J
= 0$. Hence $J = u - u \circ f$ $\mu$-a.e. for some $u
\in L^2(\mu)$ by theorem C. We obtain therefore:
\begin{equation}\label{coc}
  u - u \circ f^n (x) = \sum_{i=0}^{n-1} J \circ f^i(x) =  \log \Jac f^n(x) - \log d^{kn} .
\end{equation}
Let $\epsi>0$ and $m \geq 1$ such that $\MM:= \{\abs{u} \leq
\log m  \}$ satisfies $\mu(\MM) \geq (1-\epsi)^{1/2}$. Since $\mu$ is
mixing, $\MM_n:= \MM \cap f^{-n} \MM$ satisfies $\mu (\MM_n)
\geq \mu(\MM)^2 - \epsi \geq 1-2\epsi$ for $n$ larger than some $N'$. Let $\rho$ be small and $\tau$ be large enough such that $\mu(
\BB_n(\rho) \cap \RR_n(\tau)) \geq 1-2\epsi$ for every $n \geq 1$. We define $\TT_n:=
\BB_n(\rho) \cap \RR_n(\tau) \cap \MM_n$ and $\TT:= \limsup_{n \geq
  N'} \TT_n$. Observe that $\mu(\TT) \geq 1-4\epsi$. Let $x \in \TT$ and $(n_j)_j$ (depending on $x$) such that $x \in
\TT_{n_j}$ for every $j \geq 1$. Since $x \in \TT_{n_j} \subset
\BB_{n_j}(\rho)$, the map $f^{n_j} \circ (x + (d_x f^{n_j})^{-1})  : B(\rho)
\to \Pj^k$ is injective. \\

Let $\Lambda_{n} = d^{-n/2} \cdot
\Id_{\C^k}$. It is enough to prove that $d_x f^{n_j} = (U_j \, P_j \, V_j) \, \Lambda_{n_j}^{-1}$, where $(U_j,V_j)$ are
unitary matrices and $P_j$ is a diagonal matrix with entries in
$[a,b] \subset \R^*_+$ ($(a,b)$ being independent of $j$). Indeed, this
implies that $f^{n_j} \circ (x + \Lambda_{n_j})$ is injective on $B(\rho/b)$, completing the proof of theorem \ref{laa}. We shall omit the subscript $j$ for simplification, and denote by $\delta_1
\leq \ldots \leq \delta_k$ the singular values of $d_x f^n$. Let
$(U,V)$ be unitary matrices such that $d_x f^n = U \,
\Diag(\delta_1,\ldots,\delta_k) \, V$ (see subsection \ref{abvc}). The fact that $x \in \RR_{n}(\tau)$ yields:
\begin{equation}\label{mmm}
\delta_1  = \norm {(d_x f^n)^{-1} }^{-1} \geq d^{n/2}/\tau.
\end{equation}
Now we give an upper bound for $\delta_k$. Since $x \in \TT_n \subset \MM_n$, we have $(x,f^n(x)) \in \MM = \{ \abs{u} \leq \log m \}$. This implies by (\ref{coc}):
\[ d^{kn/2} / m \, \leq \,  \prod_{i=1}^k \delta_i = \Jac f^n(x) ^{1/2} \, \leq \, d^{kn/2} m. \]
We deduce from (\ref{mmm}):
\[ \delta_k  \, \leq  \, {\delta_1 \ldots \delta_{k-1} \over \delta_1^{k-1}} \delta_k  \, = \, {\Jac f^n(x) ^{1/2} \over  \delta_1^{k-1}} \, \leq  \,  {   d^{kn/2} m \over ( d^{n/2} / \tau )^{k-1}  }  \, =  \,  d^{n/2} \tau^{k-1}  m. \]
Thus $\Diag(\delta_1,\ldots,\delta_k) =  \Lambda_n^{-1} \, P$, where
$P$ is diagonal with entries in $[1/\tau , \tau^{k-1} m]$. We obtain
finally $d_x f^n = U \, \Lambda_n^{-1} \, P \, V =  (U \, P \, V) \,
\Lambda_n^{-1}$, as desired. \\

\emph{Acknowledgements:} We thank the referees for their careful
readings and useful suggestions. Their comments enabled us to
improve the exposition of the article.

\vspace{.5 cm}

\noindent {\footnotesize Christophe Dupont}

\vspace{.1 cm}

\noindent {\footnotesize Universit\'e Paris XI}\\
{\footnotesize CNRS UMR 8628}\\
{\footnotesize Math\'ematique, B\^atiment 425}\\
{\footnotesize 91405 Orsay Cedex, France}

\vspace{.1 cm}

\noindent{\footnotesize christophe.dupont@math.u-psud.fr}\\

\end{document}